# Corrected confidence intervals for secondary parameters following sequential tests

R. C. Weng[1] and D. S. Coad[2]

*National Chengchi University and Queen Mary, University of London*

**Abstract:** Corrected confidence intervals are developed for the mean of the second component of a bivariate normal process when the first component is being monitored sequentially. This is accomplished by constructing a first approximation to a pivotal quantity, and then using very weak expansions to determine the correction terms. The asymptotic sampling distribution of the renormalised pivotal quantity is established in both the case where the covariance matrix is known and when it is unknown. The resulting approximations have a simple form and the results of a simulation study of two well-known sequential tests show that they are very accurate. The practical usefulness of the approach is illustrated by a real example of bivariate data. Detailed proofs of the main results are provided.

## 1. Introduction

Suppose that a sequential test is carried out to compare two treatments. Then, following the test, there is interest in making valid inferences about the different parameters. For example, the primary parameter will typically be some measure of the treatment difference and there may be several secondary parameters too. These could be the individual treatment effects or the effects of baseline covariates, such as age, gender, disease stage, and so on. However, the use of a sequential design leads to the usual maximum likelihood estimators being biased and associated confidence intervals having incorrect coverage probabilities. One approach to the estimation problem is to study the approximate bias and sampling distributions of the maximum likelihood estimators.

Until recently, much of the research on estimation following sequential tests focussed on primary parameters. For example, an approach based on approximately pivotal quantities was developed by Woodroofe [24] in the context of a single stream of normally distributed observations. Here, interest lay in providing an approximate confidence interval for a mean. The work in the present paper extends this approach in several respects. We consider bivariate normal data, where interest lies in estimating the mean of the second component of the process when the first is being monitored sequentially. Further, we consider the case of an unknown covariance matrix for the process.

[1]Department of Statistics, National Chengchi University, Taipei, 11623, Taiwan, e-mail: chweng@nccu.edu.tw
[2]School of Mathematical Sciences, Queen Mary, University of London, London E1 4NS, UK, e-mail: d.s.coad@qmul.ac.uk

One of the first papers to address the problem of estimation of secondary parameters following a sequential test was [19]. For large samples, he showed how the bias of the estimator of the secondary parameter is related to that of the primary parameter, and then showed how a bias-adjusted estimator of the secondary parameter could be constructed. Gorfine has shown [7] how a theorem of Yakir [27] may be used to define an unbiased estimator of the secondary parameter. Related work has been carried out by Liu and Hall [11]. More recently, Hall and Yakir [9] develop tests and confidence procedures in a very general context.

Several authors have developed methods for the construction of confidence intervals based on approximately pivotal quantities. Whitehead, Todd and Hall show [21] how approximate confidence intervals may be obtained for a bivariate normal process when the covariance matrix is known and then show how these may be applied to problems in which approximate bivariate normality can be assumed. Liu considers [12] a similar problem and shows how the appropriate corrections may be obtained using moment expansions, though the method developed appears to be somewhat restricted. In the present paper, we consider both the known and the unknown covariance matrix cases.

The approximately pivotal quantities are constructed by considering the bivariate version of the signed root transformation, and then using a version of Stein's [15] identity and very weak expansions to determine the correction terms. The results in the known covariance matrix case are obtained by applying those of Weng and Woodroofe [17] for the two-parameter exponential family. In the unknown covariance matrix case, similar ideas to those used by Weng and Woodroofe [18] in the context of stationary autoregressive processes are used to establish the asymptotic sampling distribution of the renormalised pivotal quantity. The resulting correction terms have a simple form and complement the results of Whitehead [19].

In Section 2, the correction terms for the usual approximately pivotal quantity in the known covariance matrix case are determined using results for the two-parameter exponential family and it is shown how these may be used to construct corrected confidence intervals for the secondary parameter. The case of an unknown covariance matrix is then considered in Section 3, and the asymptotic sampling distribution of the renormalised pivotal quantity is obtained. The results of a simulation study of two well-known sequential tests are reported in Section 4 and a real example of bivariate data is used to illustrate the approach in Section 5. Some remarks and an indication of possible extensions to the present work are given in Section 6. Appendices contain detailed proofs of the main results.

## 2. Bivariate normal process with known covariance matrix

### 2.1. The general method for two-parameter exponential families

Let $X_j = (X_{1j}, X_{2j})'$ for $j = 1, \ldots, n$ be independent random vectors distributed according to a two-parameter exponential family with probability density

$$p_\theta(x) = e^{\theta' x - b(\theta)},$$

where $\theta = (\theta_1, \theta_2)' \in \Omega$ and $\Omega$ is the natural parameter space, assumed to be open. Let $L_n(\theta)$ denote the log-likelihood function based on $x_1, \ldots, x_n$, and consider the bivariate version of the signed root transformation (*e.g.* [4]) given by

$$(1) \qquad Z_{n1} = Z_{n1}(\theta) = \sqrt{2\{L_n(\hat{\theta}_n) - L_n(\tilde{\theta}_n)\}} \operatorname{sign}(\theta_1 - \hat{\theta}_{n1})$$

and

(2) $$Z_{n2} = Z_{n2}(\theta) = \sqrt{2\{L_n(\tilde{\theta}_n) - L_n(\theta)\}} \text{sign}(\theta_2 - \tilde{\theta}_{n2}),$$

where $\hat{\theta}_n = (\hat{\theta}_{n1}, \hat{\theta}_{n2})'$ is the maximum likelihood estimator and $\tilde{\theta}_n = (\theta_1, \tilde{\theta}_{n2})'$ is the restricted maximum likelihood estimator for fixed $\theta_1$. Then we have $L_n(\theta) = L_n(\hat{\theta}_n) - ||Z_n||^2/2$, where $Z_n = (Z_{n1}, Z_{n2})'$.

Consider a Bayesian model in which $\theta$ has a prior density $\xi$ with compact support in $\Omega$. Let $E_\xi$ denote expectation in the Bayesian model in which $\theta$ is replaced with a random vector $\Theta$ and let $E_\xi^n$ denote conditional expectation given $\{X_j, j = 1, \ldots, n\}$. Then the posterior density of $\Theta$ given $X_1, \ldots, X_n$ is $\xi_n(\theta) \propto e^{L_n(\theta)} \xi(\theta)$, and the posterior density of $Z_n$ is

(3) $$\zeta_n(z) \propto J(\hat{\theta}_n, \theta) \xi_n(\theta) \propto J(\hat{\theta}_n, \theta) \xi(\theta) e^{-\frac{1}{2}||z||^2},$$

where $z$ and $\theta$ are related by (1) and (2), and $J$ is a Jacobian term. From (3),

(4) $$\zeta_n(z) = f_n(z) \phi_2(z), \quad z \in \Re^2,$$

where $\phi_2$ denotes the standard bivariate normal density and

$$f_n(z) \propto J(\hat{\theta}_n, \theta) \xi(\theta).$$

Now, for $p \geq 0$, let $H_p$ be the set of all measurable functions $h : \Re^2 \to \Re$ for which $|h(z)|/c \leq 1 + ||z||^p$ for some $c > 0$, and let $H = \cup_{p \geq 0} H_p$. So, $H_0$ denotes the set of bounded functions. Let $\Phi^k$ denote the standard $k$-variate normal distribution for $k = 1, 2$ and write

$$\Gamma h = \int h d\Gamma$$

for an arbitrary measure $\Gamma$. Given $h \in H_p$, let $h_0 = \Phi^2 h$, $h_2 = h$ and

$$h_1(y_1) = \int_\Re h(y_1, w) \Phi^1(dw),$$

and

(5) $$\begin{aligned} g_1(y_1, y_2) &= e^{\frac{1}{2} y_1^2} \int_{y_1}^\infty \{h_1(w) - h_0\} e^{-\frac{1}{2} w^2} dw, \\ g_2(y_1, y_2) &= e^{\frac{1}{2} y_2^2} \int_{y_2}^\infty \{h_2(y_1, w) - h_1(y_1)\} e^{-\frac{1}{2} w^2} dw \end{aligned}$$

for $-\infty < y_1, y_2 < \infty$. Then let $Uh = (g_1, g_2)'$ and $Vh = (U^2 h + U^2 h')/2$, where $U^2 h$ is the $2 \times 2$ matrix whose $i$th column is $Ug_i$ and $g_i$ is as in (5). Lemma 1 below follows from Lemma 1 of Weng and Woodroofe [17].

**Lemma 1** (Stein's identity). *Let $r$ be a nonnegative integer. Suppose that $d\Gamma = f d\Phi^2$, where $f$ is twice differentiable on $\Re^2$ for which*

$$\int_{\Re^2} |f| d\Phi^2 + \int_{\Re^2} (1 + ||z||^r) ||\nabla f(z)|| \Phi^2(dz) < \infty$$

*and*

$$\int_{\Re^2} (1 + ||z||^r) ||\nabla^2 f(z)|| \Phi^2(dz) < \infty.$$

*Then*

$$\Gamma h = \Gamma 1 \cdot \Phi^2 h + \Phi^2(Uh)' \int_{\Re^2} \nabla f(z)\Phi^2(dz) + \int_{\Re^2} \text{tr}\{(Vh)\nabla^2 f\}d\Phi^2$$

*for all $h \in H_r$.*

From (4), the posterior distributions of $Z_n$ are of a form appropriate for Stein's identity. Let

$$\Gamma_1^\xi(\hat{\theta}_n, \theta) = \sqrt{n}\frac{\nabla f_n(Z_n)}{f_n(Z_n)}$$

and

$$\Gamma_2^\xi(\hat{\theta}_n, \theta) = n\frac{\nabla^2 f_n(Z_n)}{f_n(Z_n)}.$$

Now let $B_n$ denote the event $\{\hat{\theta}_n \in \nabla b(\Omega)\}$ and let $\Xi_0$ denote the collection of all prior densities $\xi = \xi(\theta)$ with compact support in $\Omega$ for which $\xi$ is twice differentiable almost everywhere under $P_\xi$, and $\nabla^2 \xi$ is bounded on its support. Proposition 2 below follows from Proposition 2 of Weng and Woodroofe [17].

**Proposition 2.** *Suppose that $\xi \in \Xi_0$. Then*

$$E_\xi^n\{h(Z_n)\} = \Phi^2 h + \frac{1}{\sqrt{n}}(\Phi^2 Uh)' E_\xi^n\{\Gamma_1^\xi(\hat{\theta}_n, \theta)\} + \frac{1}{n}\text{tr}[E_\xi^n\{Vh(Z_n)\Gamma_2^\xi(\hat{\theta}_n, \theta)\}]$$

*almost everywhere on $B_n$, for all $h \in H$.*

Let $N = N_a$ be a family of stopping times, depending on a design parameter $a \geq 1$. Suppose that

$$\frac{a}{N_a} \to \rho^2(\theta)$$

in $P_\theta$-probability for almost every $\theta \in \Omega$, where $\rho$ is a continuous function on $\Omega$. Suppose also that, for every compact set $K \subseteq \Omega$, there is an $\eta > 0$ such that

(6) $$P_\theta(N_a \leq \eta a) = o(a^{-q}),$$

uniformly with respect to $\theta \in K$ as $a \to \infty$, for some $q \geq 1/2$. Lemma 3 below follows from Theorem 12 of Weng and Woodroofe [17]. Moreover, by their Lemma 5 and (6) above, we have $P_\theta(B_N^c) = o(1/a)$.

**Lemma 3.** *The random vector $Z_N = (Z_{N1}, Z_{N2})'$ is uniformly integrable with respect to $P_\xi$.*

In what follows, suppose that $\theta_1$ is the primary parameter and that $\theta_2$ is a nuisance parameter. Then, for inference about $\theta_1$, it is appropriate to use $Z_{N1}$. Now, from Proposition 2,

$$E_\xi^N\{h(Z_{N1})\} = \Phi^1 h + \frac{1}{\sqrt{N}}(\Phi^1 Uh)E_\xi^N\{\Gamma_{1,1}^\xi(\hat{\theta}_N, \theta)\} + \frac{1}{N}E_\xi^N\{Vh(Z_{N1})\Gamma_{2,11}^\xi(\hat{\theta}_N, \theta)\}.$$

To determine the mean correction term for $Z_{N1}$, we take $h(z) = z$, in which case $\Phi^1 h = 0$, $\Phi^1 Uh = 1$ and $Vh(z) = 0$. Similarly, for the variance correction term, we take $h(z) = z^2$, in which case $\Phi^1 h = 1$, $\Phi Uh = 0$ and $Vh(z) = 1$. Denote by $b_{ij}$ the partial derivatives $b_{ij}(\theta) = \partial^{i+j}b(\theta)/\partial\theta_1^i\partial\theta_2^j$, and similarly for $\xi_{ij}$. Let $i_1(\theta) =$ $(b_{20} - b_{11}^2/b_{02})(\theta)$, $i_2(\theta) = b_{02}(\theta)$, $\Gamma_{1,1}^\xi(\theta, \theta) = \lim_{\omega \to \theta} \Gamma_{1,1}^\xi(\omega, \theta)$ and $\Gamma_{2,11}^\xi(\theta, \theta) = \lim_{\omega \to \theta} \Gamma_{2,11}^\xi(\omega, \theta)$, and let $\kappa(\theta)$ and $m(\theta)$ be such that

$$E_\xi\{\rho(\theta)\Gamma_{1,1}^\xi(\theta,\theta)\} = \int\int_\Omega \xi(\theta)\kappa(\theta)d\theta_1 d\theta_2 \tag{7}$$

and

$$E_\xi\{\rho^2(\theta)\Gamma_{2,11}^\xi(\theta,\theta) - 2\rho(\theta)\kappa(\theta)\Gamma_{1,1}^\xi(\theta,\theta) + \kappa^2(\theta)\} = \int\int_\Omega m(\theta)\xi(\theta)d\theta_1 d\theta_2. \tag{8}$$

Then some algebra yields

$$\kappa(\theta) = \frac{(-b_{02}, b_{11}) \cdot \nabla\rho}{b_{02} i_1^{1/2}}(\theta) \tag{9}$$
$$+ \rho(\theta)\left\{\frac{(b_{02}, -b_{11}) \cdot \nabla i_1}{6 b_{02} i_1^{3/2}}(\theta) + \frac{(b_{02}, -b_{11}) \cdot \nabla i_2}{2 b_{02}^2 i_1^{1/2}}(\theta)\right\}.$$

A similar, but more complicated expression, may also be obtained for $m(\theta)$.

Now, define

$$Z_N^{(0)} = \frac{Z_{N1} - \hat{\mu}_N^{(0)}}{\hat{\tau}_N^{(0)}}, \tag{10}$$

where

$$\hat{\mu}_N^{(0)} = \begin{cases} \hat{\kappa}_N/\sqrt{a} & \text{if } |\hat{\kappa}_N| \leq a^{1/6}\{\log(a)\}^{-1}, \\ a^{-1/3}\{\log(a)\}^{-1} & \text{if } \hat{\kappa}_N > a^{1/6}\{\log(a)\}^{-1}, \\ -a^{-1/3}\{\log(a)\}^{-1} & \text{if } \hat{\kappa}_N < -a^{1/6}\{\log(a)\}^{-1}, \end{cases} \tag{11}$$

and

$$\hat{\tau}_N^{(0)} = \begin{cases} \sqrt{1 + \hat{m}_N/a} & \text{if } |\hat{m}_N| \leq \sqrt{a}/\log(a), \\ 1 & \text{otherwise}, \end{cases} \tag{12}$$

with $\hat{\kappa}_N = \kappa(\hat{\theta}_N)$ and $\hat{m}_N = m(\hat{\theta}_N)$.

**Theorem 4.** *Let $h$ be a bounded function. Suppose that $\rho(\theta)$ is almost differentiable with respect to $\theta_1$ and $\theta_2$. If (6) holds with $q = 1$ and $\xi \in \Xi_0$, then*

$$E_\xi\{h(Z_N^{(0)})\} = \Phi^1 h + o(1/a).$$

The proof is in Appendix A.3. Theorem 4 shows that under mild conditions $Z_N^{(0)}$ is approximately standard normal to order $o(1/a)$ in the very weak sense of Woodroofe [23]. It extends Theorem 14 of Weng and Woodroofe [17] by not requiring $h$ to be symmetric and not assuming $\nabla^2 \xi$ to be continuous.

So, an asymptotic level $1 - \alpha$ confidence interval for $\theta_1$ is

$$\mathcal{I}_N = \{\theta_1 : |Z_N^{(0)}| \leq z_{\alpha/2}\}, \tag{13}$$

where $z_{\alpha/2}$ is the $100(\alpha/2)$-th percentile of the standard normal distribution.

## 2.2. The bivariate normal model with known covariance matrix

Suppose that $X_j = (X_{1j}, X_{2j})'$ for $j = 1, \ldots, n$ are independent random vectors from a bivariate normal distribution with mean vector $\theta = (\theta_1, \theta_2)'$ and covariance matrix

$$\Sigma = \begin{pmatrix} \sigma_1^2 & \gamma\sigma_1\sigma_2 \\ \gamma\sigma_1\sigma_2 & \sigma_2^2 \end{pmatrix}.$$

Let $\psi = (\sigma_1^2, \sigma_2^2, \gamma)'$. As before, let $N = N_a$ be the stopping time depending on $a$. Then, since the likelihood function is not affected by the use of a stopping time (*e.g.* [3]), the density of $X_N$ is

$$\begin{aligned}
p(x; \theta, \psi) = \exp\Bigg[&-N\log(2\pi) - \frac{N}{2}\log\{\sigma_1^2\sigma_2^2(1-\gamma^2)\} \\
&- \frac{1}{2\sigma_1^2\sigma_2^2(1-\gamma^2)}\Bigg\{\sigma_2^2\sum_{j=1}^N (x_{1j}-\theta_1)^2 + \sigma_1^2\sum_{j=1}^N (x_{2j}-\theta_2)^2 \\
&\qquad\qquad - 2\gamma\sigma_1\sigma_2\sum_{j=1}^N (x_{1j}-\theta_1)(x_{2j}-\theta_2)\Bigg\}\Bigg].
\end{aligned} \tag{14}$$

If we assume that $\theta$ is unknown and $\psi$ is known, then this model is a two-parameter exponential family with density that satisfies

$$\log p(x;\theta) = c(x) + N\theta_1 t_1 + N\theta_2 t_2 - Nb(\theta),$$

where $t_1 = \bar{x}_1/\{\sigma_1^2(1-\gamma^2)\} - \gamma\bar{x}_2/\{\sigma_1\sigma_2(1-\gamma^2)\}$, $t_2 = \bar{x}_2/\{\sigma_2^2(1-\gamma^2)\} - \gamma\bar{x}_1/\{\sigma_1\sigma_2(1-\gamma^2)\}$ and $b(\theta) = \theta'\Sigma^{-1}\theta/2$. Since $b(\theta)$ is quadratic in $\theta$, both $i_1(\theta)$ and $i_2(\theta)$ defined in Section 2.1 are constants; and therefore $\kappa(\theta)$ in (9) reduces to

$$\kappa(\theta) = \frac{(-b_{02}, b_{11}) \cdot \nabla\rho(\theta)}{b_{02} i_1^{1/2}} = -\sigma_1\rho_{10}, \tag{15}$$

where $\rho_{ij} = \partial^{i+j}\rho/\partial\theta_1^i\partial\theta_2^j$ and the second equality in (15) follows since the stopping time $N$ is assumed to depend only on $X_{11}, \ldots, X_{1N}$. Simple calculations show that the maximum likelihood estimator of $\theta$ is $(\hat{\theta}_1, \hat{\theta}_2) = (\bar{X}_{N1}, \bar{X}_{N2})$ and that the restricted maximum likelihood estimator of $\theta_2$ given $\theta_1$ is $\tilde{\theta}_2 = \tilde{\theta}_2(\theta_1) = \hat{\theta}_2 - \gamma\sigma_2(\theta_1 - \hat{\theta}_1)/\sigma_1$. By (1) and (2), it is straightforward to obtain

$$(Z_{N1}, Z_{N2}) = (\sqrt{N}\sigma_1^{-1}(\theta_1 - \hat{\theta}_1), \sqrt{N}\sigma_2^{-1}(1-\gamma^2)^{-1/2}\{\theta_2 - \hat{\theta}_2 - \gamma\sigma_2(\theta_1 - \hat{\theta}_1)/\sigma_1\}).$$

Furthermore, since the stopping time depends only on the first population, it can be shown that $m(\theta)$ in (8) satisfies

$$m(\theta) = \kappa^2(\theta) = (\sigma_1\rho_{10})^2.$$

Then, substituting these $Z_{N1}$, $\kappa$ and $m$ into (10), (11), and (12), by Theorem 4, the approximate level $1 - \alpha$ confidence interval for $\theta_1$ is as in (13).

For inference about the secondary parameter $\theta_2$, it is not appropriate to use $Z_{N2}$ as it depends on both $\theta_1$ and $\theta_2$. So, we consider the transformation

$$Z_{N1} = Z_{N1}(\theta) = \sqrt{2\{L_N(\hat{\theta}_{N1}, \hat{\theta}_{N2}) - L_N(\tilde{\theta}_{N1}, \theta_2)\}}\,\text{sign}(\theta_2 - \hat{\theta}_{N2}), \tag{16}$$

where $\tilde{\theta}_{N1} = \tilde{\theta}_{N1}(\theta_2)$ is the restricted maximum likelihood estimator of $\theta_1$ given $\theta_2$. Then $Z_{N1} = \sqrt{N}(\theta_2 - \hat{\theta}_2)/\sigma_2$. To obtain the mean correction term, we need to replace $b_{ij}$ and $\rho_{ij}$ in (15) with $b_{ji}$ and $\rho_{ji}$. So,

$$(17) \quad E_\theta(Z_{N1}) \simeq \frac{1}{\sqrt{a}}\kappa(\theta) = \frac{1}{\sqrt{a}} \frac{(-b_{20}, b_{11}) \cdot \begin{pmatrix} \rho_{01} \\ \rho_{10} \end{pmatrix}}{b_{20}(b_{02} - \frac{b_{11}^2}{b_{20}})^{1/2}}(\theta) = -\frac{1}{\sqrt{a}}\sigma_1\gamma\rho_{10}.$$

Using a similar trick, we obtain

$$(18) \quad m(\theta) = \kappa^2(\theta) = (\sigma_1\gamma\rho_{10})^2.$$

With this $Z_{N1}$ and its corresponding mean and variance corrections, we obtain a renormalised pivot $Z_N^{(0)}$ as in (10). Then, by Theorem 4, an asymptotic level $1 - \alpha$ confidence interval for $\theta_2$ is

$$(19) \quad \hat{\theta}_{N2} + \frac{\sigma_2}{\sqrt{N}}\hat{\mu}_N^{(0)} \pm \frac{\sigma_2}{\sqrt{N}}\hat{\tau}_N^{(0)} z_{\alpha/2}.$$

This interval is of the same form as the one obtained by Whitehead, Todd and Hall [21]. However, they use recursive numerical integration to calculate the correction terms instead of asymptotic approximations.

## 3. Extension to unknown covariance matrix case

In this section, we consider the following three cases:

C1. $\sigma_1$ and $\sigma_2$ are known, but $\gamma$ is unknown;
C2. $\sigma_1$ and $\sigma_2$ are unknown, but $\gamma$ is known;
C3. $\sigma_1$, $\sigma_2$ and $\gamma$ are all unknown.

When the parameters are unknown, we estimate them by $\hat{\sigma}_i^2 = \sum_{j=1}^N (X_{ij} - \hat{\theta}_i)^2 / (N-1)$ for $i = 1, 2$ and

$$\hat{\gamma} = \frac{\sum_{j=1}^N (X_{1j} - \hat{\theta}_1)(X_{2j} - \hat{\theta}_2)}{\sqrt{\sum_{j=1}^N (X_{1j} - \hat{\theta}_1)^2 \sum_{j=1}^N (X_{2j} - \hat{\theta}_2)^2}}.$$

As the main interest of this paper concerns inference about the secondary parameter $\theta_2$, in the rest of the paper we let $Z_{N1}$ be as in (16). So the corresponding $\kappa(\sigma_1, \gamma, \rho_{10})$ and $m(\sigma_1, \gamma, \rho_{10})$ are as in (17) and (18). For cases C1–C3, we consider $\hat{\kappa}_N^{(1)} = \kappa(\sigma_1, \hat{\gamma}, \hat{\rho}_{10})$, $\hat{\kappa}_N^{(2)} = \kappa(\hat{\sigma}_1, \gamma, \hat{\rho}_{10})$ and $\hat{\kappa}_N^{(3)} = \kappa(\hat{\sigma}_1, \hat{\gamma}, \hat{\rho}_{10})$, respectively; and correspondingly define $\hat{\mu}_N^{(k)}$ and $\hat{\tau}_N^{(k)}$ for $k = 1, 2, 3$ as in (11) and (12). Then, let

$$(20) \quad Z_N^{(1)} = \frac{Z_{N1} - \hat{\mu}_N^{(1)}}{\hat{\tau}_N^{(1)}}$$

and

$$(21) \quad Z_N^{(k)} = \frac{Z_{N1}(\hat{\sigma}_2) - \hat{\mu}_N^{(k)}}{\hat{\tau}_N^{(k)}}$$

for $k = 2, 3$, where $Z_{N1}(\hat{\sigma}_2) = \sqrt{N}(\theta_2 - \hat{\theta}_2)/\hat{\sigma}_2$. We will use $Z_N^{(k)}$ for $k = 1, 2, 3$ as pivotal quantities for cases C1, C2 and C3, respectively.

Define $\hat{\omega}_N = \hat{\sigma}_2^2/\sigma_2^2$. Then we can rewrite $Z_N^{(k)}$ for $k = 2, 3$ in (21) as

$$Z_N^{(k)} = \frac{(\frac{\sigma_2}{\hat{\sigma}_2})Z_{N1} - \hat{\mu}_N^{(k)}}{\hat{\tau}_N^{(k)}} = \frac{Z_{N1} - \hat{\mu}_N^{(k)}\hat{\omega}_N^{1/2}}{\hat{\omega}_N^{1/2}\hat{\tau}_N^{(k)}}. \tag{22}$$

In the rest of the paper, let $\Xi$ denote the collection of all prior densities $\xi(\psi, \theta) = \xi_1(\psi)\xi_2(\theta)$ with compact support in $(0, \infty)^2 \times (-1, 1) \times \Omega$ for which $\xi$ is twice differentiable almost everywhere under $P_\xi$ and $\nabla^2 \xi$ is bounded on its support.

**Theorem 5.** *Suppose that $\xi \in \Xi$ and that (6) holds with $q = 1$. Then, for $k = 2, 3$,*

$$\left| \int_{(0,\infty)^2 \times (-1,1)} \int_\Omega \left[ E_{\psi,\theta}\{h(Z_N^{(k)})\} - \Phi^1 h - \frac{1}{a}(\Phi_4 h)\rho^2(\theta) \right] \xi(\psi, \theta) d\theta d\psi \right| \tag{23}$$
$$= o(\frac{1}{a})$$

*for all bounded functions $h$.*

The definition of $\Phi_4$ and the proof are in Appendix A.4. Theorem 5 shows that $Z_N^{(k)}$ for $k = 2, 3$ are asymptotically distributed according to a $t$ distribution with $N$ degrees of freedom to order $o(1/a)$ in the very weak sense, since $\Phi^1 h + (\Phi_4 h)\rho^2(\theta)/a$ represents the first two terms in an Edgeworth-type expansion for the $t$ distribution (*e.g.* [1], Chap.2; [8], Chap.2). Hence,

$$P_{\psi,\theta}\{|Z_N^{(k)}| \leq z\} = 2G_N(z) - 1 + o(1/a) \tag{24}$$

very weakly, where $G_N$ denotes the $t$ distribution with $N$ degrees of freedom. So, an asymptotic level $1 - \alpha$ confidence interval for $\theta_2$ is

$$\hat{\theta}_{N2} + \frac{\hat{\sigma}_2}{\sqrt{N}}\hat{\mu}_N^{(k)} \pm \frac{\hat{\sigma}_2}{\sqrt{N}}\hat{\tau}_N^{(k)} c_{N,\alpha/2},$$

where $c_{N,\alpha/2}$ is the $100(\alpha/2)$-th percentile of the $t$ distribution with $N$ degrees of freedom. Note that the form of the above interval is similar to one obtained by Keener [10] using fixed $\theta$ expansions. However, his interval is only valid up to order $o(1/\sqrt{a})$ and only applicable to linear stopping boundaries.

The proof of Theorem 5 reveals that the correction term $(\Phi_4 h)\rho^2(\theta)/a$ in (23) arises from the use of $\hat{\omega}_N$. Since $\sigma_2$ is known for $Z_N^{(1)}$ in (20), this correction term vanishes in the asymptotic expansion for $Z_N^{(1)}$ and an immediate corollary to Theorem 5 is the following result.

**Corollary 6.** *Suppose that $\xi \in \Xi$ and that (6) holds with $q = 1$. Then*

$$\left| \int_{(0,\infty)^2 \times (-1,1)} \int_\Omega [E_{\psi,\theta}\{h(Z_N^{(1)})\} - \Phi^1 h]\xi(\psi, \theta) d\theta d\psi \right| = o(\frac{1}{a})$$

*for all bounded functions $h$.*

Therefore, $Z_N^{(1)}$ is asymptotically standard normal to order $o(1/a)$ in the very weak sense, and consequently

$$P_{\psi,\theta}\{|Z_N^{(1)}| \leq z\} = 2\Phi^1(z) - 1 + o(1/a)$$

very weakly. From this, one can set confidence intervals for $\theta_2$ as in (19), but with $\hat{\mu}_N^{(0)}$ and $\hat{\tau}_N^{(0)}$ replaced by $\hat{\mu}_N^{(1)}$ and $\hat{\tau}_N^{(1)}$.

## 4. Simulation results

### 4.1. General

Section 3 considers asymptotic results for $Z_N^{(k)}$ for a class of stopping times $N = N_a$ depending only on the first population. Specifically, let $q$ denote a measurable function on $\Re$ which is almost differentiable; let $\hat{\theta}_{n1} = \sum_{j=1}^{n} X_{1j}/n$ and

$$(25) \qquad N = N_a = \inf\{n \geq m_0 : nq(\hat{\theta}_{n1}) \geq a\} \wedge m,$$

where $m_0$ and $m$ denote the initial sample size and the maximum size, respectively; $m_0 = \lfloor a/\epsilon_0^2 \rfloor$ and $m = \lfloor a/\epsilon^2 \rfloor$, $\lfloor x \rfloor$ is the greatest integer less than or equal to $x$, $a \geq 1$ is a boundary parameter, $\epsilon$ a truncation parameter and $\epsilon_0$ controls the initial sample size. In this section, we assess the accuracy of the method for two simulated examples, the truncated sequential probability ratio test and the repeated significance test. The actual coverage probability and expected stopping time are assessed through simulation for $\sigma_1 = \sigma_2 = 1$ and selected values of $(\theta_1, \theta_2, \gamma)$.

### 4.2. Truncated sequential probability ratio test

The stopping time (25) with $q(y) = |y|$ is equivalent to

$$N = \inf\{n \geq m_0 : |S_{n1}| \geq a\} \wedge m,$$

where $S_{n1} = \sum_{j=1}^{n} X_{1j}$ is the partial sum from the first population. This is the truncated probability ratio test depending on three parameters, $a \geq 1$, $\epsilon_0$ and $\epsilon$. Simple calculations yield $a/N \to \rho^2$, where $\rho = \max\{\min(\epsilon_0, \sqrt{|\theta_1|}), \epsilon\}$. The parameter values are taken as $a = 10$, $\epsilon = \sqrt{0.1}$ and $\epsilon_0 = \sqrt{5.0}$. So, $m_0 = a/\epsilon_0^2 = 2$ and $m = a/\epsilon^2 = 100$. Tables 1 and 2 contain results for known $\sigma_1$ and $\sigma_2$, but unknown $\gamma$, that is, case C1, and for unknown $\sigma_1$, $\sigma_2$ and $\gamma$, that is, case C3, respectively.

In Table 1, we report the expected sample size and the lower and upper 0.05, 0.025 noncoverage probabilities for $Z_{N1}$ and $Z_N^{(1)}$. The results show that $Z_N^{(1)}$ is very accurate for all selected parameter values, but $Z_{N1}$ is negatively skewed. Table 2 compares the coverage probabilities using $Z_{N1}$ and $Z_N^{(3)}$. The coverage probabilities for $P_{\psi,\theta}(|Z_{N1}| \leq z_{\alpha/2})$ for $\alpha = 0.05$ and 0.1 are in the columns with the title '$Z_{N1}$'. The results using (24) for the pivotal quantity $Z_N^{(3)}$ are given under the title '$Z_N^{(3)} : t_N$'. As (23) suggests that $Z_N^{(3)}$ can be approximated by a $t$ distribution with $a/\hat{\rho}^2$ degrees of freedom, we give the results in the last two columns '$Z_N^{(3)} : t_{a/\hat{\rho}^2}$'. Apparently, the coverage probabilities for the naïve statistic are all significantly less than the nominal values. The results using $a/\hat{\rho}^2$ degrees of freedom are slightly better than those with $N$. The distribution of $Z_N^{(3)}$ shows no appreciable skewness.

TABLE 1
*Truncated sequential probability ratio test with known $\sigma_1$ and $\sigma_2$, but unknown $\gamma$;
replicates = 10,000 ($\pm$ means 1.96 standard deviations)*

|  |  | $Z_{N1}$ | | | | $Z_N^{(1)}$ | | | |
|---|---|---|---|---|---|---|---|---|---|
| $(\theta_1, \theta_2, \gamma)$ | $E_{\psi,\theta}(N)$ | L0.05 | U0.05 | L0.025 | U0.025 | L0.05 | U0.05 | L0.025 | U0.025 |
| (0.30, 1.00, 0.40) | 35.42 | 0.059 | 0.040 | 0.031 | 0.019 | 0.048 | 0.049 | 0.025 | 0.025 |
| (0.60, 1.00, 0.40) | 17.87 | 0.057 | 0.039 | 0.028 | 0.020 | 0.048 | 0.048 | 0.024 | 0.025 |
| (0.80, 1.00, 0.40) | 13.54 | 0.056 | 0.043 | 0.027 | 0.022 | 0.048 | 0.050 | 0.024 | 0.025 |
| (0.30, 1.00, 0.80) | 35.20 | 0.070 | 0.030 | 0.036 | 0.017 | 0.050 | 0.050 | 0.023 | 0.024 |
| (0.60, 1.00, 0.80) | 17.87 | 0.064 | 0.036 | 0.034 | 0.017 | 0.049 | 0.052 | 0.025 | 0.026 |
| (0.80, 1.00, 0.80) | 13.55 | 0.058 | 0.040 | 0.029 | 0.019 | 0.046 | 0.055 | 0.023 | 0.027 |
| $\pm$ |  | 0.004 | 0.004 | 0.003 | 0.003 | 0.004 | 0.004 | 0.003 | 0.003 |

TABLE 2
*Truncated sequential probability ratio test with unknown $\sigma_1$, $\sigma_2$ and $\gamma$; replicates = 10,000 ($\pm$ means 1.96 standard deviations)*

| $(\theta_1, \theta_2, \gamma)$ | $E_{\psi,\theta}(N)$ | $Z_{N1}$ 90% | $Z_{N1}$ 95% | $Z_N^{(3)} : t_N$ 90% | $Z_N^{(3)} : t_N$ 95% | $Z_N^{(3)} : t_{a/\hat{\rho}^2}$ 90% | $Z_N^{(3)} : t_{a/\hat{\rho}^2}$ 95% |
|---|---|---|---|---|---|---|---|
| (0.30, 1.00, 0.40) | 35.42 | 0.885 | 0.934 | 0.892 | 0.944 | 0.896 | 0.947 |
| (0.60, 1.00, 0.40) | 17.87 | 0.871 | 0.923 | 0.884 | 0.941 | 0.892 | 0.947 |
| (0.80, 1.00, 0.40) | 13.54 | 0.863 | 0.917 | 0.885 | 0.936 | 0.895 | 0.945 |
| (0.30, 1.00, 0.80) | 35.20 | 0.877 | 0.929 | 0.891 | 0.944 | 0.896 | 0.947 |
| (0.60, 1.00, 0.80) | 17.87 | 0.865 | 0.918 | 0.879 | 0.936 | 0.888 | 0.942 |
| (0.80, 1.00, 0.80) | 13.55 | 0.859 | 0.911 | 0.878 | 0.935 | 0.888 | 0.944 |
| $\pm$ | | 0.006 | 0.004 | 0.006 | 0.004 | 0.006 | 0.004 |

TABLE 3
*Repeated significance test with known $\sigma_1$ and $\sigma_2$, but unknown $\gamma$; replicates = 10,000 ($\pm$ means 1.96 standard deviations)*

| $(\theta_1, \theta_2, \gamma)$ | $E_{\psi,\theta}(N)$ | $Z_{N1}$ L0.05 | $Z_{N1}$ U0.05 | $Z_{N1}$ L0.025 | $Z_{N1}$ U0.025 | $Z_N^{(1)}$ L0.05 | $Z_N^{(1)}$ U0.05 | $Z_N^{(1)}$ L0.025 | $Z_N^{(1)}$ U0.025 |
|---|---|---|---|---|---|---|---|---|---|
| (0.30, 1.00, 0.40) | 75.18 | 0.064 | 0.045 | 0.034 | 0.023 | 0.052 | 0.047 | 0.026 | 0.024 |
| (0.60, 1.00, 0.40) | 27.53 | 0.061 | 0.037 | 0.031 | 0.018 | 0.047 | 0.045 | 0.023 | 0.023 |
| (0.80, 1.00, 0.40) | 16.16 | 0.060 | 0.041 | 0.032 | 0.019 | 0.050 | 0.052 | 0.024 | 0.025 |
| (0.30, 1.00, 0.80) | 74.88 | 0.093 | 0.047 | 0.049 | 0.024 | 0.052 | 0.047 | 0.025 | 0.024 |
| (0.60, 1.00, 0.80) | 27.26 | 0.083 | 0.029 | 0.041 | 0.014 | 0.051 | 0.044 | 0.025 | 0.023 |
| (0.80, 1.00, 0.80) | 16.20 | 0.067 | 0.030 | 0.032 | 0.016 | 0.049 | 0.047 | 0.025 | 0.024 |
| $\pm$ | | 0.004 | 0.004 | 0.003 | 0.003 | 0.004 | 0.004 | 0.003 | 0.003 |

### 4.3. Repeated significance test

The stopping time (25) with $q(y) = y^2$ is equivalent to

$$N = \inf\{n \geq m_0 : |S_{n1}| \geq \sqrt{na}\} \wedge m.$$

This is the repeated significance test depending on three parameters, $a \geq 1$, $\epsilon_0$ and $\epsilon$. It is easily seen that $a/N \to \rho^2$, where $\rho = \max\{\min(\epsilon_0, |\theta_1|), \epsilon\}$. We take $a = 10$, $\epsilon = \sqrt{0.1}$ and $\epsilon_0 = \sqrt{2.0}$. So, $m_0 = a/\epsilon_0^2 = 5$ and $m = a/\epsilon^2 = 100$. Tables 3 and 4 contain results for known $\sigma_1$ and $\sigma_2$, but unknown $\gamma$, that is, case C1, and for unknown $\sigma_1$, $\sigma_2$ and $\gamma$, that is, case C3, respectively.

In Table 3, we see that $Z_{N1}$ is slightly more negatively skewed than in Table 1, but $Z_N^{(1)}$ is again very accurate for all selected parameter values. The coverage probabilities in Table 4 show that the use of $Z_{N1}$ leads to significantly lower coverage probabilities than the nominal values, but using $Z_N^{(3)}$ and a $t$ distribution with $a/\hat{\rho}^2$ degrees of freedom also works very well for this test. As before, the distribution of $Z_N^{(3)}$ shows no appreciable skewness.

### 5. A practical example

In this section, we illustrate the proposed confidence interval method using the data obtained by Bellissant et al. [2]. This study was concerned with the treatment of infants of up to eight years of age suffering from gastroesophageal reflux. The infants were randomised between metoclopramide and placebo, which they received for a two-week period. The pH level in the oesophagus was measured continuously using a flexible electrode secured above the lower oesophageal sphincter. The primary response variable was the percentage reduction in acidity, measured by the proportion of time that pH $< 4$, over the two weeks of treatment.

TABLE 4
*Repeated significance test with unknown $\sigma_1$, $\sigma_2$ and $\gamma$; replicates $= 10{,}000$ ($\pm$ means 1.96 standard deviations)*

| $(\theta_1, \theta_2, \gamma)$ | $E_{\psi,\theta}(N)$ | $Z_{N1}$ | | $Z_N^{(3)} : t_N$ | | $Z_N^{(3)} : t_{a/\hat{\rho}^2}$ | |
|---|---|---|---|---|---|---|---|
| | | 90% | 95% | 90% | 95% | 90% | 95% |
| (0.30, 1.00, 0.40) | 75.18 | 0.880 | 0.934 | 0.897 | 0.947 | 0.900 | 0.948 |
| (0.60, 1.00, 0.40) | 27.53 | 0.872 | 0.925 | 0.891 | 0.939 | 0.896 | 0.946 |
| (0.80, 1.00, 0.40) | 16.16 | 0.854 | 0.907 | 0.875 | 0.933 | 0.886 | 0.942 |
| (0.30, 1.00, 0.80) | 74.88 | 0.847 | 0.911 | 0.891 | 0.945 | 0.896 | 0.948 |
| (0.60, 1.00, 0.80) | 27.26 | 0.850 | 0.908 | 0.883 | 0.938 | 0.893 | 0.945 |
| (0.80, 1.00, 0.80) | 16.20 | 0.850 | 0.904 | 0.876 | 0.934 | 0.887 | 0.945 |
| $\pm$ | | 0.006 | 0.004 | 0.006 | 0.004 | 0.006 | 0.004 |

The above variable was taken to be normally distributed and the triangular test ([20], Chap.4) was used to monitor the study. Inspections were made after groups of about four patients and the trial was stopped after the seventh interim analysis, with the conclusion that metoclopramide is not an improvement over placebo. Although Bellissant *et al.* [2] mention various normally distributed secondary response variables of interest, only standard analyses of them are carried out. For example, uncorrected confidence intervals are given for secondary parameters of interest. Thus, it is interesting to apply the corrected confidence intervals presented in Section 3 in this case.

In order to illustrate the confidence interval method, we assume that there is a single secondary response variable, the proportion of time that pH $< 4$ on day 14, and that the patients arrive in pairs, with one patient in each pair being assigned to metoclopramide and the other to placebo. The trial data give the estimates $\hat{\theta}_1 = 0.3$, $\hat{\theta}_2 = 0.07$, $\hat{\sigma}_1 = 0.5$ and $\hat{\sigma}_2 = 0.1$. To simulate the trial, we treated these values as the true values for the parameters. Further, since the sample covariance matrix was not available, we simulated the trial when $\gamma = 0.4$ and $\gamma = 0.8$, as for the two sequential tests in Section 4. As in the original trial of Bellissant *et al.* [2], we use a one-sided triangular test to test $H_0 : \theta_1 = 0$ against $H_1 : \theta_1 > 0$ and choose the design parameters so that it has significance level 5% and 95% power for $\theta_1 = 0.5$.

Let $m_a$ denote the group size, possibly depending on $a > 0$. Then the stopping time for the above triangular test is essentially of the form

$$N = \inf\{n \geq 1 : m_a | n, \ \ S_{n1}/\hat{\sigma}_1 \geq a + bn - 0.583 \ \text{ or } \ S_{n1}/\hat{\sigma}_1 \leq -a + 3bn + 0.583\},$$

where $m_a | n$ means that $m_a$ divides $n$ and $S_{n1}$ denotes the sum of the first $n$ differences in response between metoclopramide and placebo. Values are chosen for the parameters $a > 0$ and $b > 0$ in order to satisfy the error probability requirements, and the number 0.583 is a correction for overshoot of the stopping boundaries due to the discreteness of the inspection process ([20], Chap.4). Upon termination of the test, $H_0$ is rejected if $S_{N1}/\hat{\sigma}_1 \geq a + bN - 0.583$ and accepted if $S_{N1}/\hat{\sigma}_1 \leq -a + 3bN + 0.583$. Now, the above stopping time may be rewritten as

(26) $$N = \inf\{n \geq 1 : m_a | n \ \text{ and } \ nq(\hat{\theta}_{n1}/\hat{\sigma}_1) \geq a - 0.583\},$$

TABLE 5
*Triangular test with unknown $\sigma_1$, $\sigma_2$ and $\gamma$; replicates $= 10,000$ ($\pm$ means 1.96 standard deviations)*

|  |  |  | $Z_{N1}$ | | $Z_N^{(3)} : t_N$ | | $Z_N^{(3)} : t_{a/\hat{\rho}^2}$ | |
|---|---|---|---|---|---|---|---|---|
| $(\theta_1, \theta_2, \gamma)$ | Power | $E_{\psi,\theta}(N)$ | 90% | 95% | 90% | 95% | 90% | 95% |
| (0.00, 0.07, 0.40) | 0.021 | 7.43 | 0.807 | 0.864 | 0.848 | 0.921 | 0.892 | 0.935 |
| (0.00, 0.07, 0.80) | 0.021 | 7.43 | 0.815 | 0.867 | 0.857 | 0.919 | 0.896 | 0.936 |
| (0.30, 0.07, 0.40) | 0.574 | 10.49 | 0.826 | 0.885 | 0.866 | 0.927 | 0.894 | 0.949 |
| (0.30, 0.07, 0.80) | 0.574 | 10.49 | 0.780 | 0.849 | 0.860 | 0.921 | 0.892 | 0.956 |
| (0.50, 0.07, 0.40) | 0.956 | 8.17 | 0.818 | 0.877 | 0.860 | 0.926 | 0.893 | 0.942 |
| (0.50, 0.07, 0.80) | 0.956 | 8.17 | 0.812 | 0.867 | 0.859 | 0.923 | 0.896 | 0.945 |
| $\pm$ |  |  | 0.006 | 0.004 | 0.006 | 0.004 | 0.006 | 0.004 |

where $q(y) = \max(y - b, 3b - y)$. Note that (26) is a special case of more general stopping times studied by, for example, Morgan [13]. So we have $a/N \to \rho^2$, where $\rho = \max(\sqrt{\theta_1/\sigma_1 - b}, \sqrt{3b - \theta_1/\sigma_1})$, provided that $m_a = o(a)$. As in Bellissant et al. [2], we take $a = 5.495$ and $b = 0.2726$. These values may be obtained using PEST 4 [5]. Since the data are being monitored after groups of four patients, we have $m_a = 2$.

In Table 5, we report the probabilities of rejecting $H_0$, that is, the power, the expected numbers of pairs of patients, and the coverage probabilities using $Z_{N1}$ and $Z_N^{(3)}$, all of the results being based on 10,000 replications. Although the simulated sequential test satisfies the power requirement for $\theta_1 = 0.5$, it is a little conservative. This is because the above stopping time is not exactly the same as the original. Now, we know from Section 4 that the confidence intervals based on $Z_{N1}$ have coverage probabilities below the nominal values and that those based on $Z_N^{(3)}$ have roughly the correct coverage probabilities. The results in Table 5 show that the use of $Z_N^{(3)}$ leads to coverage probabilities which are usually quite close to the nominal values, especially given the small sample sizes. Note that, since our theory has been developed for the case where $\rho = \rho(\theta_1)$, when calculating the correction terms, $\sigma_1$ has been replaced with its estimate except in terms involving $\hat{\rho}$, when its true value is used. We return to this point in Section 6.

Returning to the actual trial, a standard analysis gives an uncorrected confidence interval for $\theta_2$ of (0.018, 0.122), whereas the corrected confidence interval is (0.008, 0.124) when $\gamma = 0.4$ and (0.002, 0.122) when $\gamma = 0.8$. So the approach is useful in practice, especially if the correlation coefficient is large.

## 6. Discussion

In this paper, we have shown how corrected confidence intervals for secondary parameters may be constructed following a sequential test in which one component of a bivariate normal process is being monitored. The intervals have a simple form and very weak expansions are used to justify them. Simulation of two well-known sequential tests show that the approximations are very accurate. We have also illustrated the approach using a real-life example.

We have only considered sequential tests based on the mean of the first component of a bivariate normal process. As we have seen in Section 5, a sequential test may also depend on the variance of the first component, so that $\rho = \rho(\theta_1, \sigma_1)$. The derivation of the variance correction term is more complicated in this case, since the sampling variation in $\hat{\sigma}_1^2$ needs to be allowed for. For some related work in this direction, see [26].

There may be several primary response variables in practice. So a natural extension would be to generalise the ideas in Sections 2 and 3 to a $p$-variate normal process, where $p > 2$. Such a development would require consideration of the multivariate version of the signed root transformation and an application of the results of Weng and Woodroofe [17] for the $p$-parameter exponential family in order to determine the analogues of the mean and variance correction terms in (17) and (18).

Although we have considered both the known and unknown covariance matrix cases in this paper, one assumption that we have made is that the correlation coefficient between the two components of the response vector is constant over time. This is called the proportionality case by Hall and Yakir [9]. Another natural extension would be to generalise the ideas in Sections 2 and 3 to the non-proportional case where the correlation coefficient is a function of time.

A further possible extension is to consider two binary streams of data, where the primary parameter is the log odds ratio and the secondary parameters are the individual success probabilities. Although approximations may be obtained using the results of Weng and Woodroofe [17], they do not lead to simple formulae. However, it would be interesting to compare this approach with that of Todd and Whitehead [16], and also to consider unequal sample sizes.

## Appendix A

### A.1. Wald-type equations for bivariate normal models

In this subsection, we provide some results on randomly stopped sums for the bivariate normal models. Recall the definitions of $\hat{\gamma}$ and $\hat{\sigma}_i^2$, $i = 1, 2$, in Section 3. Now define $\tilde{\sigma}_i^2 = (N-1)\hat{\sigma}_i^2/N$, $i = 1, 2$. So,

$$\tilde{\sigma}_i^2 - \sigma_i^2 = \frac{\sum_{j=1}^{N}(X_{ij} - \theta_i)^2}{N} - \sigma_i^2 - (\hat{\theta}_i - \theta_i)^2, \tag{27}$$

and $\hat{\gamma}$ defined in Section 3 can be rewritten as

$$\hat{\gamma} = \frac{\sum_{j=1}^{N}(X_{1j} - \hat{\theta}_1)(X_{2j} - \hat{\theta}_2)}{N\tilde{\sigma}_1 \tilde{\sigma}_2}. \tag{28}$$

Let $L_{1N}$ denote the likelihood function based on the first population and let $L'_{1N}$ denote the partial derivative of $L_{1N}$ with respect to $\sigma_1^2$, so that

$$L_{1N} = \exp\left\{-\frac{N}{2}\log(2\pi) - \frac{N}{2}\log(\sigma_1^2) - \frac{1}{2\sigma_1^2}\sum_{j=1}^{N}(x_{1j} - \theta_1)^2\right\}$$

and

$$L'_{1N} \equiv \frac{\partial}{\partial \sigma_1^2} L_{1N} = \left\{\frac{1}{2\sigma_1^4}\sum_{j=1}^{N}(x_{1j} - \theta_1)^2 - \frac{N}{2\sigma_1^2}\right\} L_{1N}. \tag{29}$$

Then we also have

$$\frac{\frac{\partial}{\partial \sigma_1^2}(2\sigma_1^4 L'_{1N})}{L_{1N}} = -N + \frac{1}{2\sigma_1^4}\left\{\sum_{j=1}^{N}(x_{1j} - \theta_1)^2 - \sigma_1^2 N\right\}^2. \tag{30}$$

Let $p$ be as in (14). Some of the derivations in Lemma 7 below rely on the identity (*e.g.* [22], Chap.1)

$$\text{(31)} \qquad \frac{\partial}{\partial \sigma_i^2} E_{\psi,\theta}(M_N) = \int M_N \left( \frac{\partial}{\partial \sigma_i^2} \log p \right) dP_{\psi,\theta}$$

for $i = 1, 2$, where

$$\frac{\partial}{\partial \sigma_1^2} \log p = -\frac{N}{2\sigma_1^2} + \frac{1}{2\sigma_1^2 \sigma_2^2 (1-\gamma^2)} \left\{ \frac{\sigma_2^2}{\sigma_1^2} \sum_{j=1}^N (x_{1j} - \theta_1)^2 \right.$$

$$\left. - \gamma \frac{\sigma_2}{\sigma_1} \sum_{j=1}^N (x_{1j} - \theta_1)(x_{2j} - \theta_2) \right\},$$

and $\partial \log p / \partial \sigma_2^2$ has a similar form.

**Lemma 7.** *Suppose that $\xi \in \Xi$, $M_N = M_N(X_{11}, \ldots, X_{1N}, X_{21}, \ldots, X_{2N})$, $b(\psi, \theta)$ is twice differentiable and $\nabla^2 b$ is bounded. Then the following hold:*

(i) $E_\xi \{b(\psi,\theta) M_N\} = E_\xi \{b(\psi,\theta)\} E_{\tilde{\xi}}(M_N)$, where $\tilde{\xi} = \xi b / E_\xi(b) \in \Xi$;
(ii) $E_\xi \{\sum_{j=1}^N (X_{ij} - \theta_i)^2 / N\} = E_\xi(\sigma_i^2)$ for $i = 1, 2$;
(iii) $E_\xi [\{\sum_{j=1}^N (X_{ij} - \theta_i)^2\}^2 / N] = E_\xi(2\sigma_i^4 + \sigma_i^4 N)$ for $i = 1, 2$;
(iv) $E_\xi \{\sum_{j=1}^N (X_{1j} - \theta_1)^2 \sum_{j=1}^N (X_{2j} - \theta_2)^2 / N\} = E_\xi(2\gamma^2 \sigma_1^2 \sigma_2^2 + \sigma_1^2 \sigma_2^2 N)$.

*Proof.* The proof of (i) is straightforward and hence omitted. Consider (ii). Taking $M_N = 1/N$ in (31) for $i = 1, 2$ leads to

$$\gamma E_{\psi,\theta} \left\{ \frac{\sum_{j=1}^N (X_{1j} - \theta_1)(X_{2j} - \theta_2)}{N} \right\} = -\sigma_1 \sigma_2 (1 - \gamma^2)$$

$$+ \frac{\sigma_2}{\sigma_1} E_{\psi,\theta} \left\{ \frac{\sum_{j=1}^N (X_{1j} - \theta_1)^2}{N} \right\}$$

$$- 2\sigma_1^3 \sigma_2 (1 - \gamma^2) \frac{\partial}{\partial \sigma_1^2} E_{\psi,\theta}(\frac{1}{N})$$

and a similar equation with $\sigma_1$ and $\sigma_2$ switched, where $\partial E_{\psi,\theta}(1/N)/\partial \sigma_i^2 = 0$, since $N$ defined in (25) does not depend on $\sigma_i$; and therefore,

$$\text{(32)} \qquad E_{\psi,\theta} \left\{ \frac{\sum_{j=1}^N (X_{1j} - \theta_1)^2}{N} \right\} = \frac{\sigma_1^2}{\sigma_2^2} E_{\psi,\theta} \left\{ \frac{\sum_{j=1}^N (X_{2j} - \theta_2)^2}{N} \right\}.$$

By (29),

$$E_\xi \left\{ \frac{\sum_{j=1}^N (X_{1j} - \theta_1)^2}{N} \right\} = E_\xi \left( 2\sigma_1^4 \frac{L'_{1N}}{NL_{1N}} + \sigma_1^2 \right),$$

where

$$E_\xi\left(2\sigma_1^4 \frac{L'_{1N}}{NL_{1N}}\right) = \int\int 2\sigma_1^4 \xi(\sigma_1^2,\theta_1) E_{\sigma_1,\theta_1}\left(\frac{L'_{1N}}{NL_{1N}}\right) d\sigma_1^2 d\theta_1$$

$$= -\int\int\left\{4\sigma_1^2 \xi(\sigma_1^2,\theta_1) + 2\sigma_1^2 \frac{\partial}{\partial\sigma_1^2}\xi(\sigma_1^2,\theta_1)\right\} E_{\sigma_1,\theta_1}(\frac{1}{N}) d\sigma_1^2 d\theta_1$$

$$= -\int\int 4\sigma_1^2 \xi(\sigma_1^2,\theta_1) E_{\sigma_1,\theta_1}(\frac{1}{N}) d\sigma_1^2 d\theta_1$$

$$+ \int\int 2\xi(\sigma_1^2,\theta_1)\frac{\partial}{\partial\sigma_1^2}\left\{\sigma_1^2 E_{\sigma_1,\theta_1}(\frac{1}{N})\right\} d\sigma_1^2 d\theta_1$$

$$= \int\int 2\sigma_1^4 \xi(\sigma_1^2,\theta_1)\frac{\partial}{\partial\sigma_1^2}\left\{E_{\sigma_1,\theta_1}(\frac{1}{N})\right\} d\sigma_1^2 d\theta_1,$$

where the third line follows by an integration by parts and the fact that $\xi$ is defined on a compact set and vanishes on the boundaries, and the last line is zero since $E_{\sigma_1,\theta_1}(1/N)$ does not depend on $\sigma_1$. So, $E_\xi\{\sum_{j=1}^N(X_{1j}-\theta_1)^2/N\} = E_\xi(\sigma_1^2)$, and then by (32), we have $E_\xi\{\sum_{j=1}^N(X_{2j}-\theta_2)^2/N\} = E_\xi(\sigma_2^2)$.

Now consider (iii) and (iv). First, by (30),

$$\frac{\{\sum_{j=1}^N(X_{1j}-\theta_1)^2\}^2}{N} = \frac{8\sigma_1^6 L'_{1N} + 4\sigma_1^8 L''_{1N}}{NL_{1N}} + 2\sigma_1^4 + 2\sigma_1^2\sum_{j=1}^N(X_{1j}-\theta_1)^2 - N\sigma_1^4,$$

where $E_\xi\{(8\sigma_1^6 L'_{1N} + 4\sigma_1^8 L''_{1N})/(NL_{1N})\} = 0$ can be proved in the same way as in the preceding paragraph, and $E_\xi\{\sigma_1^2\sum_{j=1}^N(X_{1j}-\theta_1)^2\} = E_\xi(\sigma_1^4 N)$ by Wald's equation. Hence,

$$(33) \quad E_\xi\left[\frac{1}{N}\left\{\sum_{j=1}^N(X_{1j}-\theta_1)^2\right\}^2\right] = E_\xi(2\sigma_1^4 + \sigma_1^4 N).$$

Next, taking $M_N = \sum_{j=1}^N(X_{1j}-\theta_1)^2/N$ in (31) for $i=1,2$ leads to

$$(34) \quad \begin{aligned}&\frac{\partial}{\partial\sigma_1^2}E_{\psi,\theta}\left\{\frac{\sum_{j=1}^N(X_{1j}-\theta_1)^2}{N}\right\}\\ &= -\frac{1}{2\sigma_1^2}E_{\psi,\theta}\left\{\sum_{j=1}^N(X_{1j}-\theta_1)^2\right\}\\ &\quad + \frac{1}{2\sigma_1^4(1-\gamma^2)}E_{\psi,\theta}\left[\frac{\{\sum_{j=1}^N(X_{1j}-\theta_1)^2\}^2}{N}\right]\\ &\quad - \frac{\gamma}{2\sigma_1^3\sigma_2(1-\gamma^2)}E_{\psi,\theta}\left[\frac{\sum_{j=1}^N(X_{1j}-\theta_1)^2\sum_{j=1}^N\{(X_{1j}-\theta_1)(X_{2j}-\theta_2)\}}{N}\right]\end{aligned}$$

and a similar equation. Taking $M_N = \sum_{j=1}^N(X_{2j}-\theta_2)^2/N$ in (31) for $i=1,2$ leads to two further equations. By (ii), (33), (34) and the other three equations, we obtain (iii) and (iv). □

*A.2. Biases and variances of estimators*

We now give some properties of the estimators for the bivariate normal models. Similar calculations are carried out by Coad and Woodroofe [6] for an adaptive normal linear model.

**Lemma 8.** *Suppose that $\xi \in \Xi$ and that (6) holds with $q = 1$. Then the following hold:*

(i) $E_\xi(\hat{\theta}_i - \theta_i) = o(1/\sqrt{a})$ *for* $i = 1, 2$;
(ii) $E_\xi(\hat{\sigma}_i^2 - \sigma_i^2) = o(1/a)$ *for* $i = 1, 2$.

*Proof.* Consider (i). First, by (6), $E_\xi[(\hat{\theta}_i - \theta_i)1_{\{a/N_a \geq 1/\eta\}}] = o(1/a)$. Then, observe that

$$\sqrt{a}(\hat{\theta}_i - \theta_i)1_{\{a/N_a < 1/\eta\}} = \sqrt{\frac{a}{N}} N^{1/2}(\hat{\theta}_i - \theta_i)1_{\{a/N_a < 1/\eta\}} \Rightarrow N(0, \sigma_i^2 \rho^2).$$

Since $a/N_a < 1/\eta$ and $N^{1/2}(\hat{\theta}_i - \theta_i)$ is uniformly integrable by Lemma 3, the left side converges in the mean. Hence (i) follows.

Next, consider (ii). From (27) and Lemma 7(ii), $aE_\xi(\tilde{\sigma}_i^2 - \sigma_i^2) = -aE_\xi\{(\hat{\theta}_i - \theta_i)^2\}$. Then, since $\hat{\sigma}_i^2 - \sigma_i^2 = \tilde{\sigma}_i^2 - \sigma_i^2 + \tilde{\sigma}_i^2/(N-1)$, we have

$$aE_\xi(\hat{\sigma}_i^2 - \sigma_i^2) = -aE_\xi\{(\hat{\theta}_i - \theta_i)^2\} + E_\xi(\frac{a}{N}\tilde{\sigma}_i^2) + o(1),$$

where $a(\hat{\theta}_i - \theta_i)^2 \Rightarrow \rho^2 \sigma_i^2 \chi_1^2$ and is uniformly integrable, by a similar argument to the preceding paragraph, and hence converges in the mean, and the second term is $E_\xi(\rho^2 \sigma_i^2) + o(1)$. So the result follows. □

A simple consequence of Lemma 8(ii) is

(35) $$E_\xi(\hat{\omega}_N - 1) = o(1/a).$$

The derivation of Lemma 9(iii)(iv) below relies on the identity

(36) $$\frac{\partial}{\partial \gamma} E_{\psi,\theta}(M_N) = \int M_N \left( \frac{\partial}{\partial \gamma} \log p \right) dP_{\psi,\theta},$$

where $p$ is as in (14) and

$$\frac{\partial}{\partial \gamma} \log p = \frac{N\gamma}{1-\gamma^2} + \frac{(1+\gamma^2)}{(1-\gamma^2)^2} \frac{\sum_{j=1}^{N}(X_{1j} - \theta_1)(X_{2j} - \theta_2)}{\sigma_1 \sigma_2}$$

$$- \frac{\gamma}{(1-\gamma^2)^2} \left\{ \frac{\sum_{j=1}^{N}(X_{1j} - \theta_1)^2}{\sigma_1^2} + \frac{\sum_{j=1}^{N}(X_{2j} - \theta_2)^2}{\sigma_2^2} \right\}$$

$$= \frac{N\gamma}{1-\gamma^2} + \frac{(1+\gamma^2)}{(1-\gamma^2)^2} \left\{ \frac{N\hat{\gamma}\hat{\sigma}_1\hat{\sigma}_2}{\sigma_1\sigma_2} + \frac{N(\hat{\theta}_1 - \theta_1)(\hat{\theta}_2 - \theta_2)}{\sigma_1\sigma_2} \right\}$$

$$- \frac{\gamma}{(1-\gamma^2)^2} \left\{ \frac{N\hat{\sigma}_1^2 + N(\hat{\theta}_1 - \theta_1)^2}{\sigma_1^2} + \frac{N\hat{\sigma}_2^2 + N(\hat{\theta}_2 - \theta_2)^2}{\sigma_2^2} \right\}.$$

**Lemma 9.** *Suppose that $\xi \in \Xi$. Then the following hold:*

(i) $E_\xi\{N(\tilde{\sigma}_i^2 - \sigma_i^2)^2\} = 2E_\xi(\sigma_i^4) + o(1)$ *for* $i = 1, 2$;
(ii) $E_\xi\{N(\tilde{\sigma}_1^2 - \sigma_1^2)(\tilde{\sigma}_2^2 - \sigma_2^2)\} = 2E_\xi(\gamma^2 \sigma_1^2 \sigma_2^2) + o(1)$;
(iii) $E_\xi(\hat{\gamma} - \gamma) = -E_\xi\{\gamma(1 - \gamma^2)/(2N)\} + o(1/a)$;
(iv) $E_\xi\{(\hat{\gamma} - \gamma)^2\} = o(1/\sqrt{a})$.

*Proof.* For (i), first from Wald's equation,

$$(37) \quad 0 = E_{\psi,\theta}\left\{\sum_{j=1}^N (X_{ij} - \theta_i)^2 - N\sigma_i^2\right\} = E_{\psi,\theta}\{N(\tilde{\sigma}_i^2 - \sigma_i^2) + N(\hat{\theta}_i - \theta_i)^2\}.$$

Next, by (27) we can write

$$N(\tilde{\sigma}_i^2 - \sigma_i^2)^2$$
$$= \frac{\{\sum_{j=1}^N (X_{ij} - \theta_i)^2\}^2}{N} - 2N(\tilde{\sigma}_i^2 - \sigma_i^2)\{\sigma_i^2 + (\hat{\theta}_i - \theta_i)^2\} - N\{\sigma_i^2 + (\hat{\theta}_i - \theta_i)^2\}^2,$$

where $E_\xi\{N(\tilde{\sigma}_i^2 - \sigma_i^2)(\hat{\theta}_i - \theta_i)^2\}$ and $E_\xi\{N(\hat{\theta}_i - \theta_i)^4\}$ are both $o(1)$. Then, together with Lemma 7(iii) and (37), we obtain the desired result.

For (ii), by (27),

$$N(\tilde{\sigma}_1^2 - \sigma_1^2)(\tilde{\sigma}_2^2 - \sigma_2^2) = \frac{1}{N}\sum_{j=1}^N (X_{1j} - \theta_1)^2 \sum_{j=1}^N (X_{2j} - \theta_2)^2 - \sigma_1^2 \sum_{j=1}^N (X_{2j} - \theta_2)^2$$
$$- \sigma_2^2 \sum_{j=1}^N (X_{1j} - \theta_1)^2 - (\hat{\theta}_1 - \theta_1)^2 \sum_{j=1}^N (X_{2j} - \theta_2)^2$$
$$- (\hat{\theta}_2 - \theta_2)^2 \sum_{j=1}^N (X_{1j} - \theta_1)^2 + N(\hat{\theta}_1 - \theta_1)^2 (\hat{\theta}_2 - \theta_2)^2$$
$$+ N\sigma_1^2 (\hat{\theta}_2 - \theta_2)^2 + N\sigma_2^2 (\hat{\theta}_1 - \theta_1)^2 + N\sigma_1^2 \sigma_2^2,$$

where we have $E_\xi\{\sum_{j=1}^N (X_{1j} - \theta_1)^2 \sum_{j=1}^N (X_{2j} - \theta_2)^2 / N\} = E_\xi(2\gamma^2 \sigma_1^2 \sigma_2^2 + \sigma_1^2 \sigma_2^2 N)$ by Lemma 7(iv), $E_\xi\{\sigma_1^2 \sum_{j=1}^N (X_{2j} - \theta_2)^2 + \sigma_2^2 \sum_{j=1}^N (X_{1j} - \theta_1)^2\} = 2E_\xi(\sigma_1^2 \sigma_2^2 N)$ by Wald's equation, $E_\xi\{N(\hat{\theta}_1 - \theta_1)^2 (\hat{\theta}_2 - \theta_2)^2\} = o(1)$ because the integrand approaches zero and is uniformly integrable, and

$$E_\xi\left[(\hat{\theta}_2 - \theta_2)^2 \left\{N\sigma_1^2 - \sum_{j=1}^N (X_{1j} - \theta_1)^2\right\}\right]$$
$$= E_\xi\left[N(\hat{\theta}_2 - \theta_2)^2 \sum_{j=1}^N \{(X_{1j} - \theta_1)^2 - \sigma_1^2\}/N\right] = o(1)$$

because $N(\hat{\theta}_2 - \theta_2)^2 = O_p(1)$ and $\sum_{j=1}^N \{(X_{1j} - \theta_1)^2 - \sigma_1^2\}/N = o_p(1)$ are both uniformly square integrable. Similarly, $E_\xi[(\hat{\theta}_1 - \theta_1)^2\{N\sigma_2^2 - \sum_{j=1}^N (X_{2j} - \theta_2)^2\}] = o(1)$. Hence (ii) follows.

Consider (iii). Taking $M_N = 1/(N\tilde{\sigma}_1\tilde{\sigma}_2)$ in (36), and then multiplying both sides by $(1 - \gamma^2)^2$ leads to

$$(1 - \gamma^2)^2 \frac{\partial}{\partial \gamma} E_{\psi,\theta} \left( \frac{1}{N\tilde{\sigma}_1 \tilde{\sigma}_2} \right) = E_{\psi,\theta} \left[ \gamma(1 - \gamma^2) \left( \frac{1}{\tilde{\sigma}_1 \tilde{\sigma}_2} \right) \right.$$

$$+ \frac{(1 + \gamma^2)}{\sigma_1 \sigma_2} \left\{ \hat{\gamma} + \frac{(\hat{\theta}_1 - \theta_1)(\hat{\theta}_2 - \theta_2)}{\tilde{\sigma}_1 \tilde{\sigma}_2} \right\}$$

(38)

$$- \frac{\gamma}{\sigma_1^2} \left\{ \frac{\tilde{\sigma}_1^2 + (\hat{\theta}_1 - \theta_1)^2}{\tilde{\sigma}_1 \tilde{\sigma}_2} \right\}$$

$$\left. - \frac{\gamma}{\sigma_2^2} \left\{ \frac{\tilde{\sigma}_2^2 + (\hat{\theta}_2 - \theta_2)^2}{\tilde{\sigma}_1 \tilde{\sigma}_2} \right\} \right].$$

From (28), the distribution of $\hat{\gamma}$ does not depend on the values of $\sigma_1$ and $\sigma_2$. So, without loss of generality, we take $\sigma_1 = \sigma_2 = 1$ in the evaluation of $E_{\psi,\theta}(\hat{\gamma})$. Letting $\psi_0 = (1, 1, \gamma)'$, then together with (38) we have

$$E_{\psi,\theta}(\hat{\gamma}) = E_{\psi_0,\theta}(\hat{\gamma})$$

$$= \frac{(1 - \gamma^2)^2}{1 + \gamma^2} \frac{\partial}{\partial \gamma} E_{\psi_0,\theta} \left( \frac{1}{N\tilde{\sigma}_1 \tilde{\sigma}_2} \right) - \frac{\gamma(1 - \gamma^2)}{1 + \gamma^2} E_{\psi_0,\theta} \left( \frac{1}{\tilde{\sigma}_1 \tilde{\sigma}_2} \right)$$

$$- E_{\psi_0,\theta} \left\{ \frac{(\hat{\theta}_1 - \theta_1)(\hat{\theta}_2 - \theta_2)}{\tilde{\sigma}_1 \tilde{\sigma}_2} \right\} + \frac{\gamma}{1 + \gamma^2} E_{\psi_0,\theta} \left( \frac{\tilde{\sigma}_1}{\tilde{\sigma}_2} + \frac{\tilde{\sigma}_2}{\tilde{\sigma}_1} \right)$$

$$+ \frac{\gamma}{1 + \gamma^2} E_{\psi_0,\theta} \left\{ \frac{(\hat{\theta}_1 - \theta_1)^2}{\tilde{\sigma}_1 \tilde{\sigma}_2} + \frac{(\hat{\theta}_2 - \theta_2)^2}{\tilde{\sigma}_1 \tilde{\sigma}_2} \right\}.$$

Now we claim that

(39) $$\frac{\partial}{\partial \gamma} E_\xi \left( \frac{1}{N\tilde{\sigma}_1 \tilde{\sigma}_2} \right) = o(\frac{1}{a}),$$

(40) $$E_\xi \left( \frac{1}{\tilde{\sigma}_1 \tilde{\sigma}_2} \right) = 1 + \frac{1}{a} \{ \frac{5}{2} E_\xi(\rho^2) + \frac{1}{4} E_\xi(\rho^2 \gamma^2 \sigma_1^2 \sigma_2^2) \} + o(\frac{1}{a}),$$

(41) $$E_\xi \left( \frac{\tilde{\sigma}_1}{\tilde{\sigma}_2} + \frac{\tilde{\sigma}_2}{\tilde{\sigma}_1} \right) = 2 + \frac{1}{a} \{ E_\xi(\rho^2) - \frac{1}{2} E_\xi(\rho^2 \gamma^2 \sigma_1^2 \sigma_2^2) \} + o(\frac{1}{a}),$$

(42) $$E_\xi \left\{ \frac{(\hat{\theta}_1 - \theta_1)^2}{\tilde{\sigma}_1 \tilde{\sigma}_2} + \frac{(\hat{\theta}_2 - \theta_2)^2}{\tilde{\sigma}_1 \tilde{\sigma}_2} \right\} = \frac{2}{a} E_\xi(\rho^2) + o(\frac{1}{a})$$

and

(43) $$E_\xi \left\{ \frac{(\hat{\theta}_1 - \theta_1)(\hat{\theta}_2 - \theta_2)}{\tilde{\sigma}_1 \tilde{\sigma}_2} \right\} = \frac{1}{a} E_\xi(\gamma \rho^2) + o(\frac{1}{a}).$$

Since the verifications of (39)–(43) are similar, here we only sketch the proof for (40). A Taylor series expansion about the point $\sigma_{10}^2 = \sigma_{20}^2 = 1$ gives

$$\frac{1}{\tilde{\sigma}_1 \tilde{\sigma}_2} \simeq 1 - \frac{1}{2}(\tilde{\sigma}_1^2 - 1) - \frac{1}{2}(\tilde{\sigma}_2^2 - 1) + \frac{3}{8}(\tilde{\sigma}_1^2 - 1)^2$$

$$+ \frac{3}{8}(\tilde{\sigma}_2^2 - 1)^2 + \frac{1}{4}(\tilde{\sigma}_1^2 - 1)(\tilde{\sigma}_2^2 - 1).$$

Then, by (i) and (ii), we obtain (40).

For (iv), we take $M_N = 1/(N\tilde{\sigma}_1\tilde{\sigma}_2)^2$ in (36), multiply both sides by $(1-\gamma^2)^2$, and then take the derivative of both sides with respect to $\gamma$. We obtain

$$\frac{\partial}{\partial \gamma}\left[(1-\gamma^2)^2 \frac{\partial}{\partial \gamma} E_{\psi,\theta}\left\{\frac{1}{(N\tilde{\sigma}_1\tilde{\sigma}_2)^2}\right\}\right]$$

$$= E_{\psi,\theta}\left(\frac{\partial}{\partial \gamma}\left[\gamma(1-\gamma^2)\left\{\frac{1}{N(\tilde{\sigma}_1\tilde{\sigma}_2)^2}\right\}\right.\right.$$

$$+ \frac{1+\gamma^2}{\sigma_1\sigma_2}\left\{\frac{\hat{\gamma}}{N\tilde{\sigma}_1\tilde{\sigma}_2} + \frac{(\hat{\theta}_1-\theta_1)(\hat{\theta}_2-\theta_2)}{N(\tilde{\sigma}_1\tilde{\sigma}_2)^2}\right\}$$

$$\left.\left.- \frac{\gamma}{\sigma_1^2}\left\{\frac{\tilde{\sigma}_1^2+(\hat{\theta}_1-\theta_1)^2}{N(\tilde{\sigma}_1\tilde{\sigma}_2)^2}\right\} - \frac{\gamma}{\sigma_2^2}\left\{\frac{\tilde{\sigma}_2^2+(\hat{\theta}_2-\theta_2)^2}{N(\tilde{\sigma}_1\tilde{\sigma}_2)^2}\right\}\right]\right)$$

$$+ E_{\psi,\theta}\left(\left[\gamma(1-\gamma^2)\left\{\frac{1}{N(\tilde{\sigma}_1\tilde{\sigma}_2)^2}\right\} + \frac{1+\gamma^2}{\sigma_1\sigma_2}\left\{\frac{\hat{\gamma}}{N\tilde{\sigma}_1\tilde{\sigma}_2} + \frac{(\hat{\theta}_1-\theta_1)(\hat{\theta}_2-\theta_2)}{N(\tilde{\sigma}_1\tilde{\sigma}_2)^2}\right\}\right.\right.$$

$$\left.\left.- \frac{\gamma}{\sigma_1^2}\left\{\frac{\tilde{\sigma}_1^2+(\hat{\theta}_1-\theta_1)^2}{N(\tilde{\sigma}_1\tilde{\sigma}_2)^2}\right\} - \frac{\gamma}{\sigma_2^2}\left\{\frac{\tilde{\sigma}_2^2+(\hat{\theta}_2-\theta_2)^2}{N(\tilde{\sigma}_1\tilde{\sigma}_2)^2}\right\}\right]\left(\frac{\partial}{\partial \gamma}\log p\right)\right)$$

$$= E_{\psi,\theta}\left\{\frac{(1+\gamma^2)^2(\hat{\gamma}-\gamma)^2}{(\sigma_1\sigma_2)^2(1-\gamma^2)^2}\right\} + I_a,$$

where $E_\xi(I_a) = O(1/a)$ and the last equality follows from tedious calculations, which are omitted here. Since

$$E_\xi\left(\frac{\partial}{\partial \gamma}\left[(1-\gamma^2)^2 \frac{\partial}{\partial \gamma} E_{\psi,\theta}\left\{\frac{1}{(N\tilde{\sigma}_1\tilde{\sigma}_2)^2}\right\}\right]\right) = o(\frac{1}{\sqrt{a}}),$$

the required result follows. □

Note that, in the absence of a stopping time, Lemma 9(iii) reduces to the usual bias formula for the sample correlation coefficient in the fixed-sample case (*e.g.* [14], Chap.5).

**Lemma 10.** *Suppose that $\xi \in \Xi$. Then $\sqrt{N}(\hat{\omega}_N - 1) \Rightarrow N(0,2)$. Moreover, $N(\hat{\omega}_N - 1)^2$ is uniformly integrable with respect to $P_\xi$.*

*Proof.* The first statement follows since

$$\sqrt{N}(\hat{\omega}_N - 1) = \frac{\sqrt{N}}{\sigma_2^2}(\hat{\sigma}_2^2 - \sigma_2^2)$$

$$= \sqrt{N}\frac{\sum_{j=1}^{N}\{(X_{2j}-\theta_2)^2 - \sigma_2^2\}}{(N-1)\sigma_2^2} + \frac{\sqrt{N}}{N-1} - \frac{\sqrt{N}(\hat{\theta}_2-\theta_2)^2}{\sigma_2^2},$$

where the first term on the right-hand side converges in distribution to $N(0,2)$ by Anscombe's theorem and the last two terms are $o_p(1)$.

For the second statement, it suffices to show that $N(\hat{\omega}_N - 1)^2$ converges to $2\chi_1^2$ in the mean. From Lemma 9(i) and the relationship between $\hat{\sigma}_2$ and $\tilde{\sigma}_2$, we have $E_\xi\{N(\hat{\omega}_N - 1)^2\} \to 2$. So the result follows. □

Two additional results are needed for Lemma 11. By (6) with $q \geq 1/2$ and Lemma 9(i), we have

(44) $\quad E_\xi[a(\tilde{\sigma}_i^2 - \sigma_i^2)^2 1_{\{N_a > \eta a\}}] = E_\xi\left[\frac{a}{N}N(\tilde{\sigma}_i^2 - \sigma_i^2)^2 1_{\{N_a > \eta a\}}\right] = O(1)$

for $i = 1, 2$, and, by (37), we have

$$(45) \quad E_\xi[a||\hat{\theta} - \theta||^2 1_{\{N_a > \eta a\}}] = E_\xi\left[\frac{a}{N}N||\hat{\theta} - \theta||^2 1_{\{N_a > \eta a\}}\right] = O(1).$$

**Lemma 11.** *Let $g(\psi, \theta)$ be twice continuously differentiable on a compact set $K \subseteq (0, \infty)^2 \times (-1, 1) \times \Omega$. Suppose that $\xi \in \Xi$ and (6) holds with $q \geq 1/2$. Then*

$$E_\xi[\{g(\psi, \theta) - g(\hat{\psi}, \hat{\theta})\}1_{\{N_a > \eta a\}}] = o(1/\sqrt{a}).$$

*Proof.* By compactness and continuity, there exists $C > 0$ such that

$$|g(\psi, \theta) - g(\hat{\psi}, \hat{\theta}) - (\psi - \hat{\psi}, \theta - \hat{\theta})' \nabla g(\hat{\psi}, \hat{\theta})| \leq C(||\hat{\psi} - \psi||^2 + ||\hat{\theta} - \theta||^2).$$

Now, since $E_\xi[(||\hat{\psi} - \psi||^2 + ||\hat{\theta} - \theta||^2)1_{\{N_a > \eta a\}}] = o(1/\sqrt{a})$ by Lemma 9(iv), (44) and (45), and $||E_\xi(\hat{\theta} - \theta)|| + ||E_\xi(\hat{\psi} - \psi)|| = o(1/\sqrt{a})$ by Lemma 8(i)(ii) and Lemma 9(iii), the statement follows by using the arguments in Proposition 6.13 of Weng and Woodroofe [18]. □

Note that, if $K = \cup_{i=1}^q K_i$, $K_i$ are compact sets, $K_i^o \cap K_j^o = \emptyset$ for $i \neq j$, where $K_i^o$ denotes the interior of $K_i$, and $g$ is twice piecewise continuously differentiable on $K_i$, then we can write

$$E_\xi\{g(\psi, \theta) - g(\hat{\psi}, \hat{\theta})\} = \sum_{i=1}^q \int_{K_i} \xi(\psi, \theta) E_{\psi,\theta}\{g(\psi, \theta) - g(\hat{\psi}, \hat{\theta})\} d\theta d\psi$$

$$= \sum_{i=1}^q \frac{1}{c_i} \int_{K_i} \xi_i(\psi, \theta) E_{\psi,\theta}\{g(\psi, \theta) - g(\hat{\psi}, \hat{\theta})\} d\theta d\psi$$

$$= \sum_{i=1}^q \frac{1}{c_i} E_{\xi_i}\{g(\psi, \theta) - g(\hat{\psi}, \hat{\theta})\},$$

where $c_i$ are normalising constants and $\xi_i = c_i \xi 1_{K_i}$. Thus, Lemma 11 holds for such $g$. In particular, it applies to $\kappa$ and yields $E_\xi[\{\hat{\kappa}_N^{(k)} - \kappa\}1_{\{N_a > \eta a\}}] = o(1/\sqrt{a})$ for $k = 0, 1, 2, 3$.

### A.3. Proof of Theorem 4

Three lemmas are required for the proof.

**Lemma 12.** *Let $h$ be a bounded function, and let*

$$H_0(\sigma, \mu) = \int_\Re h(\frac{z - \mu}{\sigma})\phi(z) dz$$

*and*

$$H_1(\sigma, \mu) = \int_\Re z h(\frac{z - \mu}{\sigma})\phi(z) dz$$

*for $\sigma > 0$ and $-\infty < \mu < \infty$. Then $H_0$ and $H_1$ have continuous derivatives of all orders. Further, at $\mu = 0$ and $\sigma = 1$, we have $H_0 = \Phi^1 h$, $\partial H_0/\partial \mu = -\Phi^1 Uh$, $\partial H_0/\partial \sigma = -2\Phi^1 Vh$, $\partial^2 H_0/\partial \mu^2 = 2\Phi^1 Vh$, $H_1 = \Phi^1 Uh$, $\partial H_1/\partial \mu = -2\Phi^1 Vh$, $\partial H_1/\partial \sigma = 0$ and $\partial^2 H_1/\partial \mu^2 = 0$.*

This lemma is a simple extension of Lemma 13 of Weng and Woodroofe [17] and can be proved analogously. Note that, if $h$ is symmetric, then $\Phi^1 U h = 0$, and hence $\partial H_0 / \partial \mu = H_1 = 0$, and Lemma 12 reduces to their Lemma 13.

Now define

$$R_{2,a}(h) = a \left( E_\xi^N \{h(Z_N)\} - \Phi^2 h - \frac{1}{\sqrt{N}} (\Phi^2 U h)' E_\xi^N \{\Gamma_1^\xi(\theta, \theta)\} \right.$$
$$\left. - \frac{1}{a} \mathrm{tr}[\Phi^2 V h E_\xi^N \{\rho^2(\theta) \Gamma_2^\xi(\theta, \theta)\}] \right).$$

Then, by Proposition 2, we have

$$R_{2,a}(h) = \sqrt{a} \left[ \sqrt{\frac{a}{N}} (\Phi^2 U h)' E_\xi^N \{\Gamma_1^\xi(\hat{\theta}, \theta) - \Gamma_1^\xi(\theta, \theta)\} \right]$$
$$+ \frac{a}{N} \mathrm{tr}[E_\xi^N \{V h(Z_N) \Gamma_2^\xi(\hat{\theta}, \theta)\}] - \mathrm{tr}[\Phi^2 V h E_\xi^N \{\rho^2(\theta) \Gamma_2^\xi(\theta, \theta)\}]$$
$$= R_{2,a}^{(1)}(h) + R_{2,a}^{(2)}(h).$$

Lemma 13 below is similar to Theorem 7 of Weng and Woodroofe [17], but here we consider $R_{2,a}(h)$ for all bounded $h$, not necessarily symmetric.

**Lemma 13.** *If (6) holds with $q = 1$ and $\xi \in \Xi_0$, then $\lim_{a \to \infty} |E_\xi[R_{2,a}(h) 1_{\{N_a > \eta a\}}]| = 0$ for all bounded $h$.*

*Proof.* First, $\lim_{a \to \infty} |E_\xi \{R_{2,a}^{(2)}(h)\}| = 0$ follows by the same argument used to prove Theorem 7 of Weng and Woodroofe [17]. Next, since

$$|E_\xi[R_{2,a}^{(1)}(h) 1_{\{N_a > \eta a\}}]|$$
$$= \sqrt{a} \left| E_\xi((\Phi^2 U h)' E_\xi^N \left[ \sqrt{\frac{a}{N}} \{\Gamma_1^\xi(\hat{\theta}, \theta) - \Gamma_1^\xi(\theta, \theta)\} \right] 1_{\{N_a > \eta a\}}) \right|$$
$$\leq C_1 \sqrt{a} E_\xi \{\|\Gamma_1^\xi(\hat{\theta}, \theta) - \Gamma_1^\xi(\theta, \theta)\|\}$$

for some constant $C_1$, the right-hand side is $o(1)$ by Lemma 11. □

**Proof of Theorem 4.** Since $h$ is bounded and both $P_\xi(N_a \leq \eta a)$ and $P_\xi(B_N^c)$ are $o(1/a)$, it suffices to show that $E_\xi[h(Z_N^{(0)}) 1_{\{N_a > \eta a\} \cap B_N}] = \Phi^1 h + o(1/a)$. Write $h(Z_N^{(0)}) = h_a(Z_{N1})$. Then, by the definition of $R_{2,a}$,

$$E_\xi^N \{h(Z_N^{(0)})\} = E_\xi^N \{h_a(Z_{N1})\}$$
(46)
$$= \Phi^1 h_a + \frac{1}{\sqrt{a}} (\Phi^1 U h_a) E_\xi^N \{\rho(\theta) \Gamma_{1,1}^\xi(\theta, \theta)\}$$
$$+ \frac{1}{a} (\Phi^1 V h_a) E_\xi^N \{\rho^2(\theta) \Gamma_{2,11}^\xi(\theta, \theta)\} + \frac{1}{a} R_{2,a}(h_a),$$

where $E_\xi[R_{2,a}(h_a) 1_{\{N_a > \eta a\}}] \to 0$ as $a \to \infty$ by Lemma 13. Since $h$ here may not be symmetric, by Lemma 12 two additional terms arise in the analysis of (46), namely,

$$A(h) = \frac{1}{\sqrt{a}} (\Phi^1 U h) E_\xi^N \{\rho(\theta) \Gamma_{1,1}^\xi(\theta, \theta)\} - \frac{1}{\sqrt{a}} (\Phi^1 U h) \hat{\kappa}_N^{(0)}.$$

To show that the effect of non-symmetry of $h$ vanishes, observe that

$$E_\xi[A(h) 1_{\{N_a > \eta a\}}] = \frac{1}{\sqrt{a}} (\Phi^1 U h) E_\xi[\{\rho(\theta) \Gamma_{1,1}^\xi(\theta, \theta) - \hat{\kappa}_N^{(0)}\} 1_{\{N_a > \eta a\}}]$$
$$= \frac{1}{\sqrt{a}} (\Phi^1 U h) E_\xi[\{\kappa(\theta) - \hat{\kappa}_N^{(0)}\} 1_{\{N_a > \eta a\}}],$$

where the last line is $o(1/a)$ by Lemma 11. So, the theorem follows. □

### A.4. Proof of Theorem 5

Given a measurable function $h$, $s > 0$, $c > 0$ and $\nu \in \Re$, let

(47) $$h^*(z) = h\{s^{-\frac{1}{2}}c^{-1}(z-\nu)\},$$

$$\Psi_0(h;\nu,s) = -(\Phi^1 Uh)\nu + (\Phi^1 Vh)\nu^2 - 2(\Phi^1 Vh)(c-1)$$
$$-(\Phi^1 Vh)(s-1) - (\Phi_3 h)\nu(s-1) + \frac{1}{2}(\Phi_4 h)(s-1)^2$$

and

$$\Psi_1(h;\nu,s) = -2(\Phi^1 Vh)\nu + (\Phi_3 h)(s-1),$$

where

$$\Phi_3 h = \frac{1}{2}\int_\Re (2-z^2)zh(z)\Phi^1\{dz\}$$

and

$$\Phi_4 h = \int_\Re \{\frac{1}{4}(z^2-1)^2 - \frac{1}{2}\}h(z)\Phi^1\{dz\}.$$

**Lemma 14.** *There is a constant $C$ for which*

$$|\Phi^1 h^* - \Phi^1 h - \Psi_0(h;\nu,s)| \le C\{|\nu|^3 + |s-1|^3 + |c-1|^{3/2}\},$$

$$|\Phi^1 Uh^* - \Phi^1 Uh - \Psi_1(h;\nu,s)| \le C\{|\nu|^2 + |s-1|^2 + |c-1|\}$$

*and*

$$|\Phi^1 Vh^* - \Phi^1 Vh| \le C\{|\nu| + |s-1| + |c-1|\},$$

*for all $|\nu| \le 1, |s-1| \le 1/2, |c-1| \le 1/2$ and bounded $h$.*

We omit the proof of this lemma since it can be derived in a similar manner to Lemma 1 of Woodroofe and Coad [25].

**Proof of Theorem 5.** We shall only consider $Z_N^{(2)}$, as the same argument applies to $Z_N^{(3)}$. From (22), we can write $h\{Z_N^{(2)}\} = h^*(Z_{N1})$, where $Z_{N1}$ is defined in (16) and $h^*(z)$ is as in (47). As in the proof of Theorem 4, we only need to consider the set $\{N_a > \eta a\} \cap B_N$. If $\psi = (\sigma_1^2, \sigma_2^2, \gamma)'$ is known, the bivariate normal model is a two-parameter exponential family. Let $E_{\xi_2}^{N,\psi}$ denote the conditional expectation given $\psi$ and the data by time $N$. So, by Proposition 2,

$$E_{\xi_2}^{N,\psi}\{h^*(Z_{N1})\} - \Phi^1 h = E_{\xi_2}^{N,\psi}\{h^*(Z_{N1})\} - \Phi^1 h^* + \Phi^1 h^* - \Phi^1 h$$
(48) $$= \frac{1}{\sqrt{a}}(\Phi^1 Uh^*)E_{\xi_2}^{N,\psi}\{\rho(\theta)\Gamma_{1,1}^\xi(\theta,\theta)\}$$
$$+ \frac{1}{a}(\Phi^1 Vh^*)E_{\xi_2}^{N,\psi}\{\rho^2(\theta)\Gamma_{2,11}^\xi(\theta,\theta)\}$$
$$+ \frac{1}{a}R_{2,a}\{h^*(Z_N)\} + \Phi^1 h^* - \Phi^1 h,$$

where $E_{\xi_2}[R_{2,a}(h^*)1_{\{N_a>\eta a\}}] \to 0$ as $a \to \infty$ by Lemma 13. Then, by Lemma 14 we can write the last two lines of (48) as

$$\frac{1}{\sqrt{a}}\{\Phi^1 Uh - 2(\Phi^1 Vh)\hat{\omega}_N^{1/2}\hat{\mu}_N^{(2)} + (\Phi_3 h)(\hat{\omega}_N - 1)\}E_{\xi_2}^{N,\psi}\{\rho(\theta)\Gamma_{1,1}^{\xi}(\theta,\theta)\}$$
$$+ \frac{1}{a}(\Phi^1 Vh)E_{\xi_2}^{N,\psi}\{\rho^2(\theta)\Gamma_{2,11}^{\xi}(\theta,\theta)\} - (\Phi^1 Uh)\hat{\omega}_N^{1/2}\hat{\mu}_N^{(2)} + (\Phi^1 Vh)\{\hat{\omega}_N^{1/2}\hat{\mu}_N^{(2)}\}^2$$
$$- 2(\Phi^1 Vh)\{\hat{\tau}_N^{(2)} - 1\} - (\Phi^1 Vh)(\hat{\omega}_N - 1) - (\Phi_3 h)\hat{\omega}_N^{1/2}\hat{\mu}_N^{(2)}(\hat{\omega}_N - 1)$$
$$+ \frac{1}{2}(\Phi_4 h)(\hat{\omega}_N - 1)^2 + o(\frac{1}{a})$$
$$= \frac{1}{\sqrt{a}}(\Phi^1 Uh)I_N + \frac{1}{a}(\Phi^1 Vh)II_N + \frac{1}{\sqrt{a}}(\Phi_3 h)III_N + (\Phi_4 h)IV_N + o(\frac{1}{a}),$$

where

$$I_N = E_{\xi_2}^{N,\psi}\{\rho(\theta)\Gamma_{1,1}^{\xi}(\theta,\theta)\} - \hat{\omega}_N^{1/2}\hat{\kappa}_N^{(2)},$$

$$II_N = E_{\xi_2}^{N,\psi}\{\rho^2(\theta)\Gamma_{2,11}^{\xi}(\theta,\theta)\} - 2\hat{\omega}_N^{1/2}\hat{\kappa}_N^{(2)}E_{\xi_2}^{N,\psi}\{\rho(\theta)\Gamma_{1,1}^{\xi}(\theta,\theta)\} + \{\hat{\omega}_N^{1/2}\hat{\kappa}_N^{(2)}\}^2$$
$$-2a\{\hat{\tau}_N^{(2)} - 1\} - a(\hat{\omega}_N - 1),$$

$$III_N = E_{\xi_2}^{N,\psi}[(\hat{\omega}_N - 1)\{\rho(\theta)\Gamma_{1,1}^{\xi}(\theta,\theta) - \hat{\omega}_N^{1/2}\hat{\kappa}_N^{(2)}\}]$$

and

$$IV_N = \frac{1}{2}(\hat{\omega}_N - 1)^2.$$

To prove (23), it suffices to show that $E_{\xi}[I_N 1_{\{N_a>\eta a\}}] = o(1/\sqrt{a})$, that $E_{\xi}[II_N 1_{\{N_a>\eta a\}}] = o(1)$, that $E_{\xi}[III_N 1_{\{N_a>\eta a\}}] = o(1/\sqrt{a})$ and that $aE_{\xi}[IV_N \times 1_{\{N_a>\eta a\}}] = E_{\xi}\{\rho^2(\theta)\} + o(1)$. For $I_N$, recall from (7) that we may write $E_{\xi_2}\{\rho(\theta) \times \Gamma_{1,1}^{\xi}(\theta,\theta)\} = E_{\xi_2}\{\kappa(\sigma_1,\gamma,\rho_{10})\}$, which together with Lemma 11 yields

$$\sqrt{a}E_{\xi}[I_N 1_{\{N_a>\eta a\}}] = \sqrt{a}E_{\xi}[\{\kappa - \hat{\omega}_N^{1/2}\hat{\kappa}_N^{(2)}\}1_{\{N_a>\eta a\}}] = o(1).$$

Next, consider $II_N$. By consistency of $\hat{\sigma}_2$ and $\hat{\kappa}_N^{(2)}$, and (18),

$$E_{\xi}[\rho^2(\theta)\Gamma_{2,11}^{\xi}(\theta,\theta) - 2\rho(\theta)\hat{\kappa}_N^{(2)}\hat{\omega}_N^{1/2}\Gamma_{1,1}^{\xi}(\theta,\theta) + \{\hat{\kappa}_N^{(2)}\hat{\omega}_N^{1/2}\}^2]$$
$$= E_{\xi}\{\rho^2(\theta)\Gamma_{2,11}^{\xi}(\theta,\theta) - 2\rho(\theta)\kappa\Gamma_{1,1}^{\xi}(\theta,\theta) + \kappa^2\} + o(1)$$
$$= E_{\xi}(\kappa^2) + o(1).$$

So, by definition of $\hat{\tau}_N^{(2)}$ and (35),

$$E_{\xi}(II_N) = E_{\xi}[\kappa^2 - 2a\{\hat{\tau}_N^{(2)} - 1\}] - aE_{\xi}(\hat{\omega}_N - 1) + o(1)$$
$$= -aE_{\xi}(\hat{\omega}_N - 1) + o(1)$$
$$= o(1).$$

For $III_N$, write

$$\sqrt{a}(III_N) = \sqrt{a}(\hat{\omega}_N - 1)\{\rho(\theta)\Gamma_{1,1}^{\xi}(\theta,\theta) - \hat{\omega}_N^{1/2}\hat{\kappa}_N^{(2)}\},$$

where $|\rho(\theta)\Gamma_{1,1}^{\xi} - \hat{\omega}_N^{1/2}\hat{\kappa}_N^{(2)}|$ is bounded, and $\sqrt{a}(\hat{\omega}_N - 1)$ converges to a limit with mean zero and is uniformly integrable by Lemma 10. So $E_{\xi}(III_N) = o(1/\sqrt{a})$.

For $IV_N$, first observe by (6) that $E_{\xi}[IV_N 1_{\{N_a\leq\eta a\}}] = o(1/a)$. Then, note that we have $aIV_N 1_{\{N_a>\eta a\}} = a(\hat{\omega}_N - 1)^2 1_{\{N_a>\eta a\}}/2$, which is uniformly integrable and approaches $\rho^2 \chi_1^2$ by Lemma 10. So, the desired result follows. □

**Acknowledgments**

Part of this work was carried out while the first author was visiting the University of Sussex during July and August 2003, and in receipt of Study Visit Grant 15600 from The Royal Society. The first author was also partially supported by the National Science Council of Taiwan. The authors are grateful to Professor J.R. Whitehead for suggesting the data used in Section 5. They also wish to thank the two referees for their comments, which have led to an improved paper.